# The Cycle Problem:
# An Intriguing Periodicity to the Zeros of the Riemann Zeta Function

David D. Baugh
dbaugh@rice.edu

ABSTRACT. Summing the values of the real portion of the logarithmic integral of $n^{\rho}$, where $\rho$ is one of a consecutive series of zeros of the Riemann zeta function, reveals an unexpected periodicity to the sum. This is the cycle problem.

As part of a larger project, the value of the logarithmic integral, $\int_0^x 1/\log t \, dt$, was evaluated for $x$ equals $n^{\rho}$, where $\rho$ is a zero of the Riemann zeta function. While looking at the summed real portion of the values of the integral with $n$ equal to $10^6$ after every million zeros, an apparent pattern developed, see Fig. 1. This is designated as the cycle problem because these repeating periods were unexpected and are not fully explained.

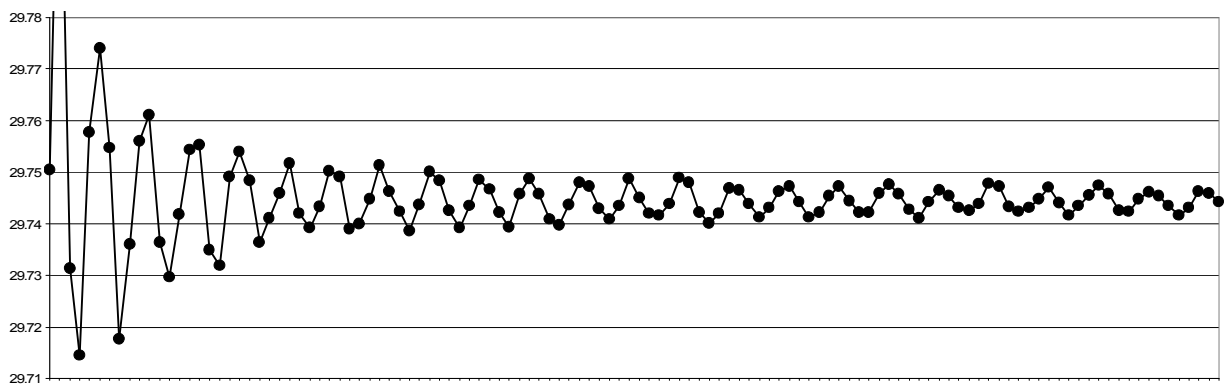

Figure 1: Plot that started this investigation

This pattern can be approximated by a fairly simple formula, see Fig. 2. It is basically sinusoidal with a decreasing amplitude and increasing period.

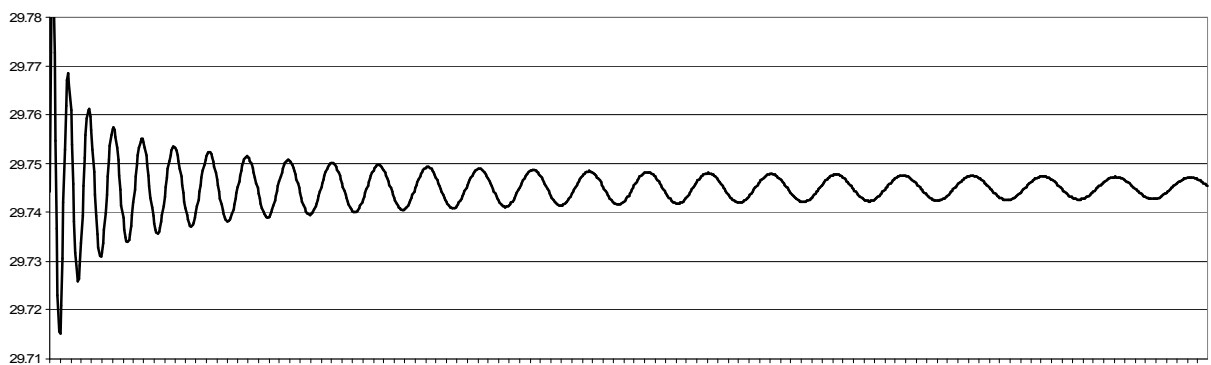

Figure 2: `Plot[29.745 + Sin[x`$^{5/8}$`]/(3 x`$^{5/8}$`),{x,6 Pi,1004 Pi]`

---

In exploring this phenomenon in detail for smaller values of *n*, it was discovered that the pattern can be much more complex than a simple sine wave. This example uses an *n* of 1295 and the first 30,000 Riemann zeta zeros, see Fig. 3.

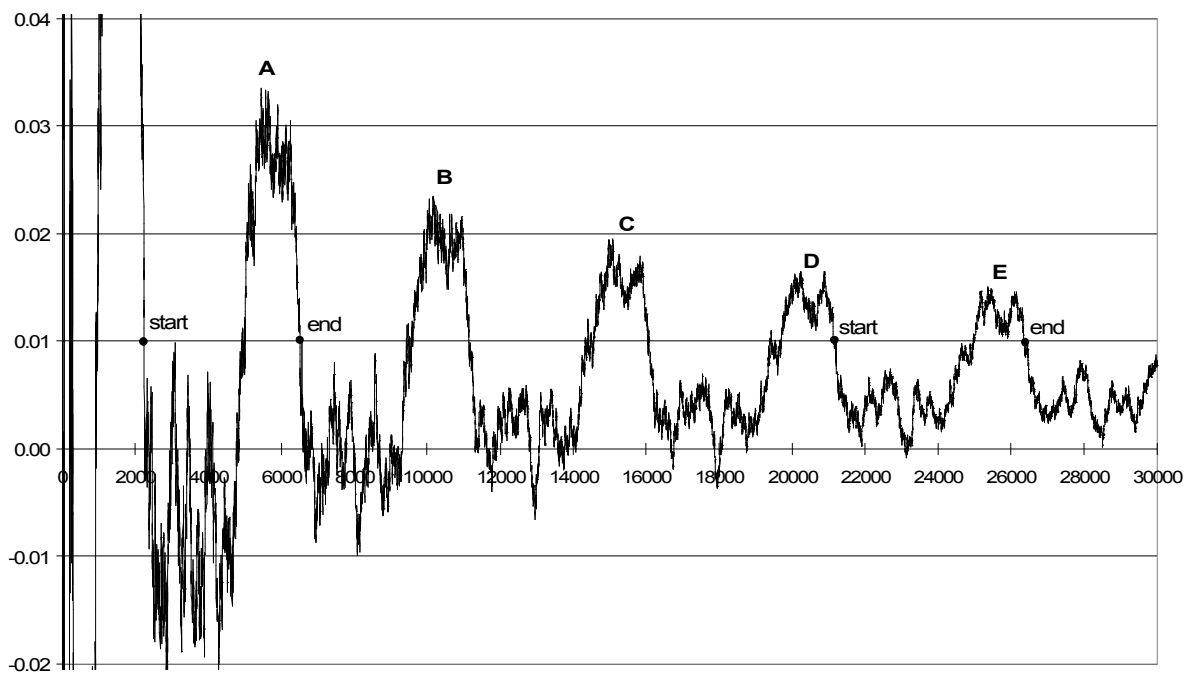

Figure 3: Accumulated real portion using 30,000 zeros and n = 1295

Not every value of *n* produces such an obvious pattern. An example of this is shown in the following figure which uses an *n* of 1302 and the first 30,000 Riemann zeta zeros, see Fig. 4.

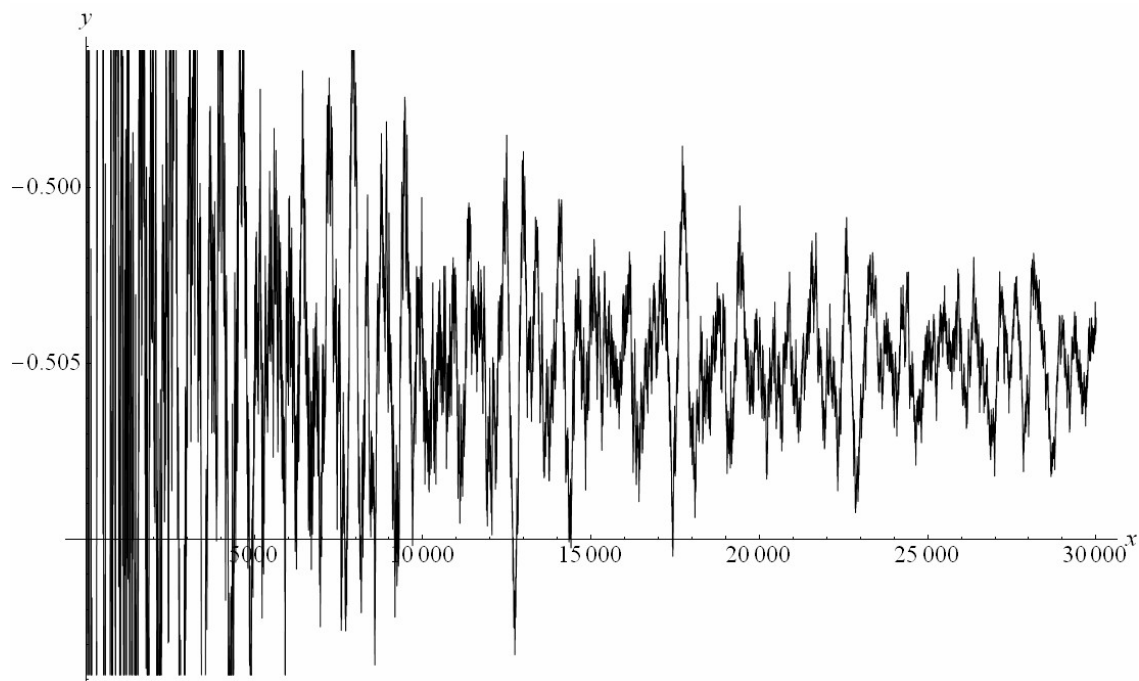

Figure 4: Accumulated real portion using 30,000 zeros and n = 1302

It should be noted that summing the imaginary portion does produce an obvious pattern for this value, see Fig. 5.

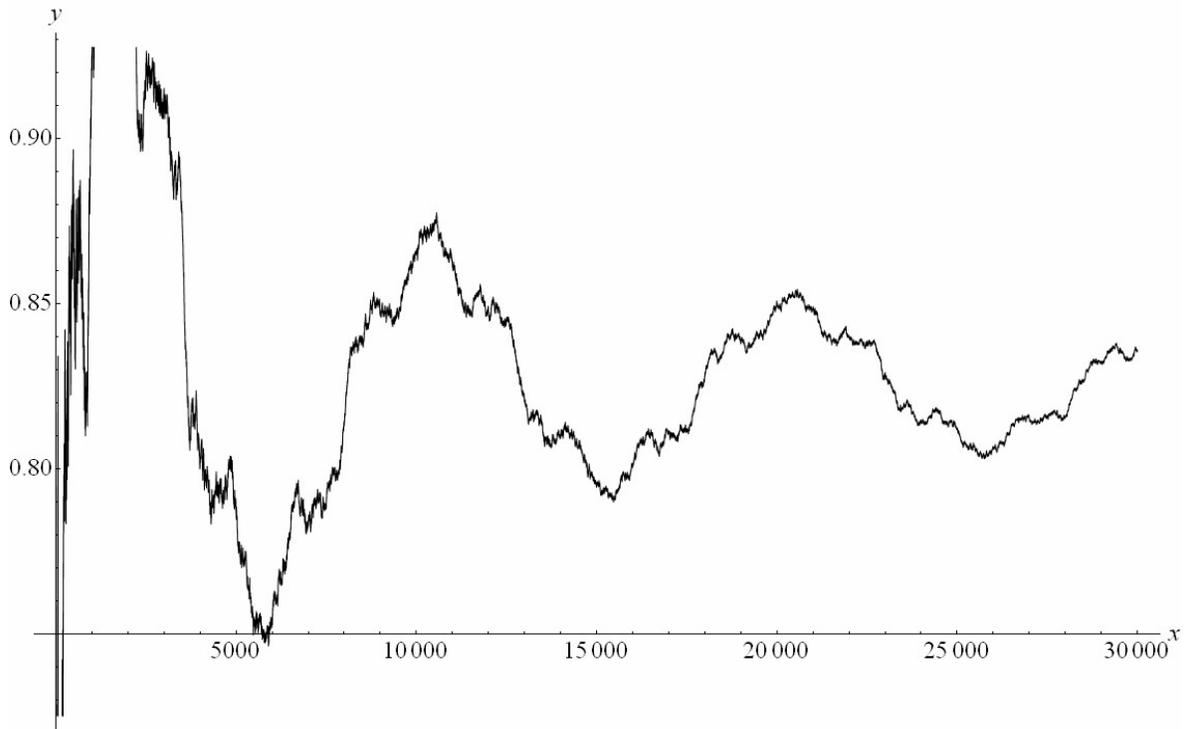

Figure 5: Accumulated imaginary portion using 30,000 zeros and n = 1302

In every case the values are bounded by a narrowing envelope. This is due to a feature of the logarithmic integral for *n* raised to a complex power. When plotted for some value of *n*, here the number 12, it produces a spiral approaching $\pi i$. The values produced using the first 20 Riemann zeta zeros are represented as large dots along this zoom of the spiral, see Fig. 6.

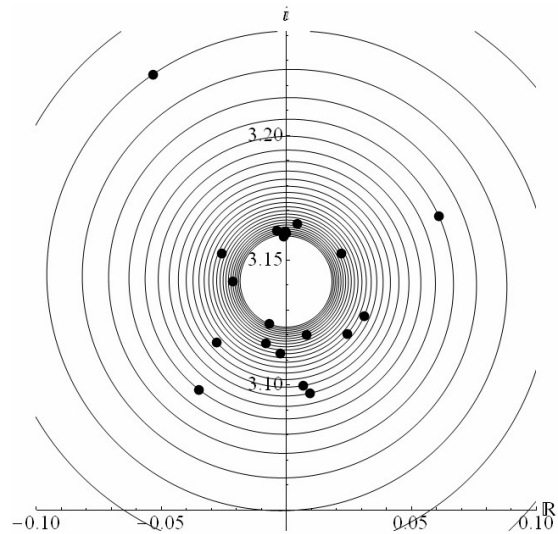

Figure 6: Partial plot of logarithmic integral

It is this spiraling inward that determines the characteristic envelope, see Fig. 7. Here each dot represents the real portion of the logarithmic integral of $12^{\rho}$ where $\rho$ is one of the first 4,000 Riemann zeta zeros. The plot does not display the full vertical scale.

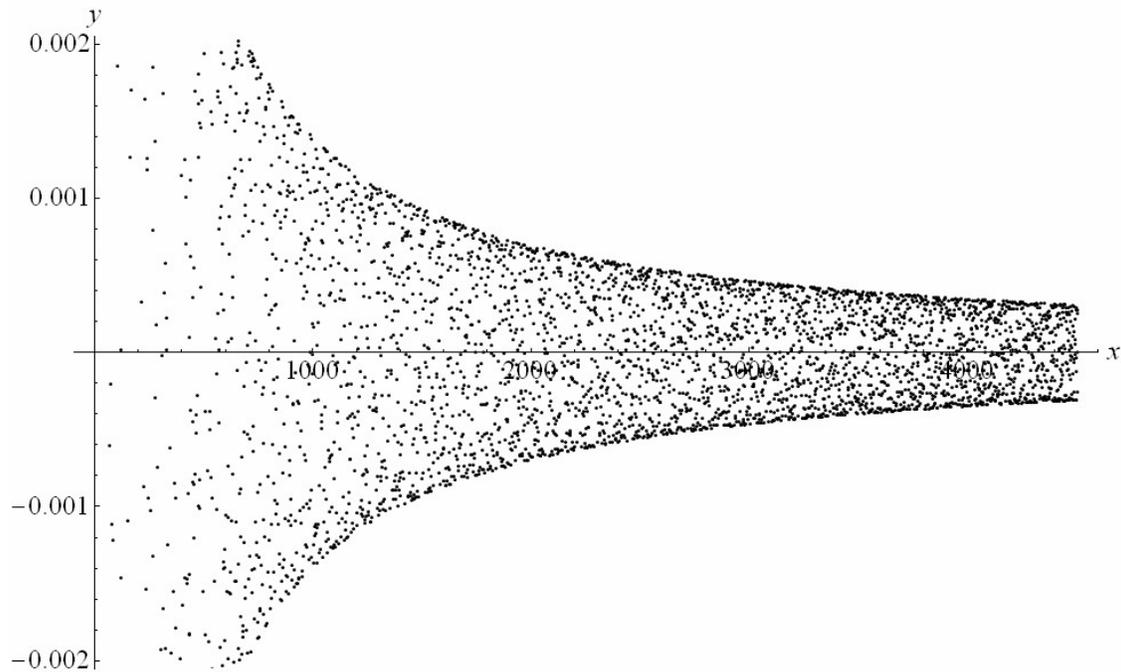

Figure 7: Real portion using 4,000 zeros and n = 12

Plotting these points as a cumulative sum generates an even more unexpected pattern, see Fig. 8.

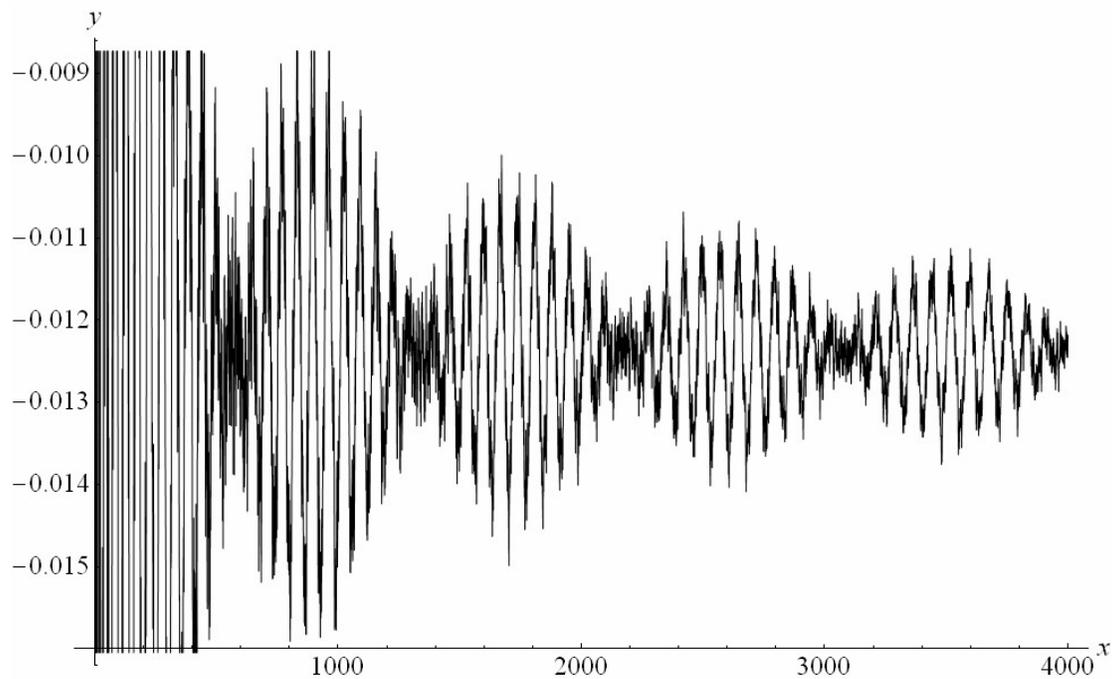

Figure 8: Accumulated real portion using 4,000 zeros and n = 12

For a pattern such as the one in Fig. 3 to occur requires a remarkable coincidental distribution of these dots given that their position on the spiral is determined by the zeros of the Riemann zeta function, see Fig. 9.

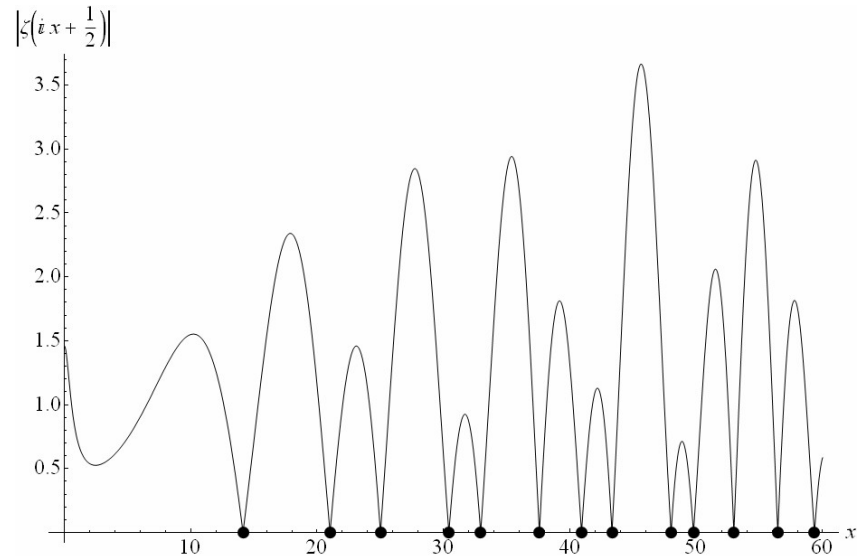

Figure 9: Riemann zeta function with zeros indicated

Period A in Fig. 3 is generated with an *n* of 1295 and the Riemann zeta zeros 2249 through 6481. Period E uses zeros 21133 through 26331. The additional 966 zeros cause the elongation of the period. The larger values of the zeros cause the attenuation. The crux of the cycle problem is to understand how these two sets of Riemann zeta zeros could produce two such similar constellations of dots on the logarithmic integral plot.

Riemann referred to the function being investigated as the "periodic terms" [1] of the Riemann prime counting function. His use of the term periodic meant only that the individual values oscillated between positive and negative. He was not implying the periodicity observed here. The cycle problem suggests that there is a readily discernable periodicity to the zeros of the Riemann zeta function. This information is presented for discussion in the hope that better minds will discover the underlying explanation for the patterns observed even if it is something that should have been obvious to the researcher.

All logarithmic integral calculations were made with *Mathematica 6* using the function `ExpIntegralEi`. Although substantial effort has been expended to assure that the periodicity observed is genuine and not just an artifact of the calculations, this remains a possibility. Many thanks go to Dr. David Damanik of Rice University for his supervision and Dr. Michael Rubinstein of the University of Waterloo for his zeros.

---

2000 *Mathematics Subject Classification*. Primary 11Y40; Secondary 11M26.